\documentclass[oneside,11pt]{amsart}
\usepackage{amssymb,latexsym,amsmath,amsthm,enumitem,mathrsfs}
\usepackage[text={390pt,615pt},centering]{geometry}
\usepackage{hyperref}
\usepackage{fancyhdr}
\usepackage{color}
\usepackage{stmaryrd}
\usepackage{mathrsfs}
\usepackage{lineno}
\usepackage{enumitem}
\usepackage{soul}
\usepackage[normalem]{ulem}

\newtheorem{thm}{Theorem}
\newtheorem*{thm*}{Theorem}
\newtheorem{prop}[thm]{Proposition}

\theoremstyle{remark}

\newcommand{\N}{\mathbb{N}}
\newcommand{\R}{\mathbb{R}}

\newcommand{\Z}{\mathbb{Z}}

\newcommand{\ep}{\varepsilon}

\newcommand{\beq}{\begin{equation}}
\newcommand{\eeq}{\end{equation}}

\makeatletter
\def\@tocline#1#2#3#4#5#6#7{\relax
  \ifnum #1>\c@tocdepth % then omit
  \else
    \par \addpenalty\@secpenalty\addvspace{#2}%
    \begingroup \hyphenpenalty\@M
    \@ifempty{#4}{%
      \@tempdima\csname r@tocindent\number#1\endcsname\relax
    }{%
      \@tempdima#4\relax
    }%
    \parindent\z@ \leftskip#3\relax \advance\leftskip\@tempdima\relax
    \rightskip\@pnumwidth plus4em \parfillskip-\@pnumwidth
    #5\leavevmode\hskip-\@tempdima
      \ifcase #1
       \or\or \hskip 1em \or \hskip 2em \else \hskip 3em \fi%
      #6\nobreak\relax
    \hfill\hbox to\@pnumwidth{\@tocpagenum{#7}}\par% <---- \dotfill -> \hfill
    \nobreak
    \endgroup
  \fi}
\makeatother

\definecolor{pink}{rgb}{1,.2,.6}
\definecolor{orange}{rgb}{0.7,0.3,0}
\definecolor{blue}{rgb}{.2,.6,.75}
\definecolor{green}{rgb}{.4,.7,.4}
\definecolor{purple}{RGB}{127,0,255}

\begin{document}

\title[A new type of superorthogonality]{A new type of superorthogonality}

\author[Gressman]{Philip T. Gressman}
\address{University of Pennsylvania, 209 South 33rd Street, Philadelphia PA 19104}
\email{gressman@math.upenn.edu}

\author[Pierce]{Lillian B. Pierce}
\address{Duke University, 120 Science Drive, Durham NC 27708 }
\email{pierce@math.duke.edu}

\author[Roos]{Joris Roos}
\address{University of Massachusetts Lowell, Lowell, MA 01854
}
\email{joris\_roos@uml.edu}

\author[Yung]{Po-Lam Yung}
\address{ Australian National University, Canberra, ACT 2601 
\& The Chinese University of Hong Kong, Shatin, Hong Kong}
\email{polam.yung@anu.edu.au, \, plyung@math.cuhk.edu.hk}

\maketitle

\begin{abstract}
We provide a simple criterion on a family of functions that implies a square function estimate on $L^p$ for every even integer $p \geq 2$. This defines a new type of superorthogonality that is verified by checking a less restrictive criterion than any other type of superorthogonality that is currently known.
\end{abstract}

\section{Statement of the results}
\subsection{Square function estimate}
Let $\{f_j\}_{j \in J}$ be a family of complex-valued functions on a measure space $(X,d\mu)$, with indices ranging over a countable set $J$.  Here and in all that follows, $L^p$ norms indicate $L^p(X,d\mu)$. 
A square function estimate of the form
\[ \| \sum_{j \in J} f_j   \|_{L^{p}} \leq C_{p} \| (\sum_{j \in J} |f_j|^2 )^\frac{1}{2} \|_{L^{p}}\]
is an essential ingredient in many methods in harmonic analysis. 
We prove a simple criterion under which  this square function estimate holds for each even integer $p \geq 2$.

\begin{thm}\label{thm_main}
Let $r \geq 1$ be an integer. Suppose a family $\{f_j\}_{j \in J}$ of   functions has the property that
\beq\label{vanish} \int f_{j_1} \cdots f_{j_r} \overline{f_{j_{r+1}} \cdots f_{j_{2r}}}  = 0 \eeq
whenever $j_1,\ldots,j_{2r}$ are all distinct. If $ (\sum_{j\in J} |f_j|^2 )^\frac{1}{2}$ belongs to $L^{2r}$, the series $\sum_{j \in J} f_j$ converges unconditionally in $L^{2r}$, and 
\beq\label{square} \| \sum_{j\in J} f_j   \|_{L^{2r}} \leq C_{r} \| (\sum_{j\in J} |f_j|^2 )^\frac{1}{2} \|_{L^{2r}}. \eeq
In particular, we may take $C_r=1$ for $r=1$ and $C_r \leq 2^{1/2}((2r)!-1)^{1/2}
$ for $r \geq 2$. 
\end{thm} 

Note that in our terminology, a family of functions  could be an ordered sequence $\{f_j\}_{j\in \N}$, but we also allow $\{f_j\}_{j \in J}$ where $J$ is an unordered countable set.
 
 Theorem \ref{thm_main} immediately implies that any family of functions satisfying the hypothesis of the theorem obeys an $L^{2r} \ell^2$ decoupling estimate,
\[  \| \sum_{j\in J} f_j\|_{L^{2r}} \leq C_r (\sum_{j\in J} \|f_j\|_{L^{2r}}^2)^{1/2}\]
for the same constant $C_r$ as in (\ref{square}), since by Minkowski's inequality for $L^p$ norms the square function norm is majorized by the right-hand side shown above.

We prove the theorem  via an application of a pointwise inequality, which may be stated in terms of sequences of complex numbers.

\begin{prop}\label{prop_numerical}
Let $k \geq 2$ be an integer and suppose $a^1,\ldots,a^k$ are finite or absolutely convergent infinite sequences of complex numbers with the terms of $a^i$ denoted $a^i_j$, where $j$ ranges over a set $J$ of indices. Let $s$ denote the sum operation,  
\beq\label{sum_operation} s(a^i) := \sum_{j\in J} a^i_j \eeq
and let
\[ Q_k(a^1,\ldots,a^k) := \sum_{\substack{j_1,\ldots,j_k\in J \\ j_1,\ldots,j_k \text{ distinct}}} a^1_{j_1} \cdots a^k_{j_k}. \]
Let
\begin{align*}
    A & := \max_{i=1,\ldots,k} |s(a^i)|, \\
    B & := \max_{i=1,\ldots,k} ||a^i||_{\ell^2}.
\end{align*}
Then
\begin{equation} \left|Q_k(a^1,\ldots,a^k) - \prod_{i=1}^k s(a^i) \right| \leq (k!-1) B^2 \left( \max\{A,B\} \right)^{k-2}. \label{theinequality} \end{equation}
\end{prop}
It is important in our application that the term $A$ is not the $\ell^1$ norm, but instead the sum operation (\ref{sum_operation}), which does not take absolute values inside the sum.
In particular, to deduce Theorem \ref{thm_main}, we will apply the consequence  that if $\prod_{i=1}^k s(a^i)$ is known to be real, then  (\ref{theinequality}) implies
 \beq\label{real_inequality}
 \prod_{i=1}^k s(a^i)
 \leq (k!-1) B^2 \left( \max\{A,B\} \right)^{k-2} + \Re(Q_k(a^1,\ldots,a^k)).
 \eeq

\subsection{Superorthogonality}
Theorem \ref{thm_main} identifies a new type of superorthogonality that   is broader than any   previously identified type, in the sense that it is verified by checking a less stringent criterion than any previously identified type. A family of functions $\{f_j\}_{j \in J}$ is said to be superorthogonal for $2r$-tuples of a certain Type if 
\beq\label{vanish'} \int f_{j_1} \cdots f_{j_r} \overline{f_{j_{r+1}} \cdots f_{j_{2r}}}  = 0 \eeq
whenever the tuple of indices $j_1,\ldots,j_{2r}$ lies in a certain subset of $J^{2r}$.  The framework of superorthogonality was developed in \cite{Pie20}, which identified several types:

\begin{itemize}
\item Type I*: (\ref{vanish'}) vanishes when $j_1,\ldots, j_r$ is not a permutation of $j_{r+1},\ldots, j_{2r}$;

\item Type I: (\ref{vanish'}) vanishes when an index $j_i$ appears an odd number of times in $j_1,\ldots,j_{2r}$; 

\item Type II: (\ref{vanish'}) vanishes when an index $j_i$ appears precisely once in $j_1,\ldots,j_{2r}$. 
\end{itemize}
 
\noindent
When the indices are positive integers, one can also define:
\begin{itemize}
    \item Type III: (\ref{vanish'}) vanishes when an index $j_i$ is strictly larger than all other indices in $j_1,\ldots,j_{2r}$.
\end{itemize}
 For example, the Rademacher functions $\{r_j\}_j$ satisfy the Type I condition for $2r$-tuples for all integers $r \geq 1$ \cite[\S2.1]{Pie20}. Any  
family $\{f_j\}_{j \in J}$ in which the members $f_j$  are mutually independent random variables, and
each has mean zero, satisfies the Type II condition for $2r$-tuples for all integers $r \geq 1$ \cite[\S 1]{Pie20}.  Let $\{w_j\}_j$ denote the sequence of Walsh-Paley functions, and let $P_jf = \sum_{2^{j-1}\leq m < 2^j}c(f)w_j$ denote a dyadic partial sum of the Walsh-Paley expansion for a given function $f$. Then the dyadic partial sums   $\{P_jf\}_j$  satisfy the Type III condition for $2r$-tuples for all integers $r \geq 1$ \cite[\S 4.3]{Pie20}.

Any family of functions that is Type I*, I, or II for $2r$-tuples satisfies a square function estimate on $L^{2r}$. See \cite{Pie20}, which also observed that  any family of functions that is Type III for $2r$-tuples and obeys two auxiliary properties (e.g. a maximal estimate, motivated by a method due to Paley) satisfies a square function estimate on $L^{2r}$.
In fact, proving that a family of Type I* functions satisfies a square function estimate (with $C_r \leq (r!)^{1/2r}$) is quite simple. As the set of indices upon which the Type assumes (\ref{vanish'}) vanishes shrinks, the abstract verification of the square function estimate becomes more tricky.

The principle of using superorthogonality to prove a square function estimate can be traced, retrospectively, through a wide variety of papers over the past 90 years;   \cite{Pie20} identified its role in: Rademacher functions and Khintchine's inequality (an application of Type I); Paley's proof \cite{Pal32} of the $L^p$-norm convergence of Walsh-Paley series  (Type III); Ionescu and Wainger's influential method \cite{IW} for bounding discrete singular Radon transforms  (Type II);   the present authors' work \cite{GGPRY20}  on square function estimates for non-degenerate curves such as the moment curve (Type I*);   results in number theory such as Burgess's celebrated bound  \cite{Bur63A}  for short character sums (Type II); and Fouvry-Kowalski-Michel's striking proof \cite{FKM15}  of square-root cancellation for sums of products of trace functions, as a consequence of the proof of the Weil conjectures (Type I).

\subsection{Type IV superorthogonality} 
The present work identifies a new type of superothogonality:
\begin{itemize}
    \item 
Type IV: (\ref{vanish'}) vanishes when all indices $j_1,\ldots,j_{2r}$ are distinct.
\end{itemize}
Theorem \ref{thm_main} proves unconditionally that any family of functions that is superorthogonal of Type IV for $2r$-tuples satisfies a square function estimate on $L^{2r}$.

It is striking that the criterion for Type IV is the least restrictive of any previously identified type. 
That is to say, if $J$ is the index set of the sequence $\{f_j\}_{j\in J}$ and $Z(\mathrm{Type}) \subset J^{2r}$ denotes the  set of tuples  for  which  (\ref{vanish'}) must vanish for that  Type to hold, then 
\[ Z(\mathrm{IV}) \subsetneq Z(\mathrm{III})  \subsetneq Z(\mathrm{II})  \subsetneq Z(\mathrm{I})  \subsetneq Z(\mathrm{I^*}).\]
Equivalently, 
\[ \text{$\{f_j\}_{j\in J}$ of Type I*} \Rightarrow \text{Type I}\Rightarrow \text{Type II}\Rightarrow \text{Type III}\Rightarrow \text{Type IV}.\]
(Of course, if the indices are not positive integers, Type III is omitted.)
In particular, as a result of Theorem \ref{thm_main}, it is now unconditionally true that any family of functions that is Type III superorthogonal satisfies the square function estimate (\ref{square}).

We exhibit several families of real-valued functions that are Type IV but not any other type.
 We first remark that there is an elementary construction of a sequence of $N$ real-valued functions $f_1,\ldots,f_N$ such that 
\[
 \int_{[0,1]} f_{j_1}f_{j_2}\cdots f_{j_N} =0
\]
 precisely when $j_1,j_2,\ldots,j_N$ are distinct.  
 To see this, divide the unit interval into $N$ subintervals $I_j=[j/N, (j+1)/N)$. Define $f_j = \chi_{^cI_j}$ where the complement is defined with respect to the unit interval, that is $^cI_j=[0,1]\setminus I_j$. Then the sequence $\{f_j\}_{j=0}^{N-1}$ of $N$ functions is Type IV  superorthogonal on $[0,1]$ for $N$-tuples. 
 
 Our next theorem enhances  this construction using Rademacher functions to create an infinite sequence of functions on $\R$ that is Type IV for $2k$-tuples.

\begin{thm}\label{thm_example}
Fix any positive integer $k$. There exist real-valued, piecewise-constant functions $\{f_i\}_{i=1}^\infty$ on $\R$ such that
\begin{itemize}
\item For each $i \in {\mathbb N}$, $f_i \in L^{2k}(\R)$.
\item If $i_1,\ldots,i_{2k} \in {\mathbb N}$ and $i_{s} \neq i_{s'}$ whenever $s \neq s'$, then
\begin{equation} \int_{\R} f_{i_1}(t) \cdots f_{i_{2k}}(t)\,dt = 0. \label{innerp0} \end{equation}
\item If $i_1,\ldots,i_{2k} \in {\mathbb N}$ and $i_{s} = i_{s'}$ for some $s \neq s'$, then
\begin{equation} \int_{\R} f_{i_1}(t) \cdots f_{i_{2k}}(t)\,dt > 0. \label{innerp1} \end{equation}
\end{itemize}
\end{thm}

Another interesting example  is given by Haar functions.
Let $\mathcal{D}$ denote the set of standard dyadic intervals $I=[2^k \ell, 2^k (\ell+1))$ in $\R$ with $k,\ell\in\mathbb{Z}$.
Each $I\in\mathcal{D}$ is associated with the Haar function 
\beq\label{Haar_dfn}
\psi_I = |I|^{-1/2} (\mathbf{1}_{I_{l}} - \mathbf{1}_{I_r}), \eeq
where $I_l, I_r\in\mathcal{D}$ denote the left and right children of $I$, respectively.

\begin{prop}\label{prop_Haar}
The family of Haar functions $\{\psi_I\}_{I \in \mathcal{D}}$, indexed by the set of dyadic intervals, is superorthogonal of Type IV on $\R$ for $2r$-tuples, for each integer $r\geq 1$. It is  not superorthogonal of Type II or III on $\R$ for $2r$-tuples, for any $r\geq 2$. However, the subfamily of Haar functions associated with dyadic intervals contained in a given compact interval is superorthogonal of Type III for $2r$-tuples, for each $r \geq 1$.
\end{prop}
 Consequently, applying Theorem \ref{thm_main} for a given integer $r \geq 1$ to Haar functions one recovers the classical square function estimate 
\begin{equation}\label{eqn:sqfcthaar}
\| f \|_{L^{p}} \lesssim \Big\| \Big( \sum_{k\in\mathbb{Z}} |\mathbb{D}_k f|^2 \Big)^{1/2} \Big\|_{L^{p}}
\end{equation}
for $p=2r$, where $\mathbb{D}_k f$ denotes the dyadic martingale differences
\[ \mathbb{D}_k f = \sum_{I\in\mathcal{D}, |I|=2^{-k}} \langle f,\psi_I\rangle \psi_I. \]
The estimate \eqref{eqn:sqfcthaar} holds for all $p\in (1,\infty)$ and is usually proved via Calder\'{o}n--Zygmund theory and the converse square function estimate (see e.g. \cite[p. 14, Thm. 1.27]{Per01}).

\section{Proof of Proposition \ref{prop_numerical}}
The pointwise inequality in Proposition \ref{prop_numerical} may be proved as follows. When $k=2$, the proposition is essentially a minor variation on Cauchy-Schwarz because
\begin{align*}
    \left| \sum_{j_1,j_2\in J, j_1 \neq j_2} a^1_{j_1} a^2_{j_2} - \sum_{j_1,j_2\in J} a^1_{j_1} a^2_{j_2} \right| & = \left| \sum_{j_1\in J} a^1_{j_1} a^2_{j_1} \right| \leq ||a^1||_{\ell^2} ||a^2||_{\ell^2} \leq B^2.
\end{align*}
For larger values of $k$, scaling implies that one may assume without loss of generality that $||a^i||_{\ell^2} \leq 1$ for each $i$ with equality for at least one $i$ (since $||a^i||_{\ell^2} = 0$ for all $i$ is a trivial case). Under this normalization assumption, it suffices to prove that
\[ \left|Q_k(a^1,\ldots,a^k) - \prod_{i=1}^k s(a^i) \right| \leq C_k  \left( \max\{1,|s(a^1)|,\ldots,|s(a^k)| \} \right)^{k-2} \]
for each $k > 2$ and appropriate constants $C_k$.
The quantities $Q_{k+1}$ and $Q_{k}$ are related by the recursive identity
\begin{align*} Q_{k+1}(a^1,\ldots,a^{k+1}) = Q_{k}(a^1,\ldots,a^{k}) s(a^{k+1}) & - Q_{k} (a^1 a^{k+1},a^2,\ldots,a^{k}) \\ & - \cdots - Q_{k}(a^1,\ldots,a^{k} a^{k+1}) \end{align*}
when $k \geq 2$.
(Here, for sequences $a=\{a_1,a_2,\ldots\}$ and $b=\{b_1,b_2,\ldots\}$ we let $ab$ denote the sequence $\{a_1b_1,a_2b_2,\ldots\}$.) By induction and the triangle inequality,
\begin{align*} |Q_k(a^1 a^{k+1},a^2,& \ldots,a^k)| \leq |s(a^1a^{k+1})| |s(a^2)| \cdots |s(a^k)| \\ & + C_k \left( \max\{1, |s(a^1 a^{k+1})|, |s(a^2)|,\ldots, |s(a^k)|\} \right)^{k-2}. \end{align*}
By the normalization condition on each $a^i$, $|s(a^1 a^{k+1})| \leq 1$, so the right-hand side is 
\[  \leq |s(a^2)| \cdots |s(a^k)|  + C_k \left( \max\{1 , |s(a^2)|,\ldots, |s(a^k)|\} \right)^{k-2}.\]
We can conclude that certainly
\[ |Q_k(a^1 a^{k+1},a^2, \ldots,a^k)| \leq (C_k + 1) \left( \max\{1, |s(a^1)|,\ldots, |s(a^{k+1})|\} \right)^{k-1} \]
and likewise for $|Q_k (a^1,a^2 a^{k+1},\ldots,a^k)|$, etc. It follows by the triangle inequality that
\begin{align*}
    \left|Q_{k+1}(a^1,\ldots,a^{k+1}) - \prod_{i=1}^{k+1} s(a^i)\right | & \leq \\ |s(a^{k+1})| & \left|Q_k(a^1,\ldots,a^k) - \prod_{i=1}^k s(a^i) \right| \\
    + k (C_k + 1)& \left( \max\{1, |s(a^1)|,\ldots, |s(a^{k+1})|\} \right)^{k-1}. 
\end{align*}
By induction, it follows for each $k \geq 2$ that
\begin{align*}  \left|Q_{k+1}(a^1,\ldots,a^{k+1}) - \prod_{i=1}^{k+1} s(a^i)\right | & \\  \leq ((k+1) C_k + k) & \left( \max\{1, |s(a^1)|,\ldots, |s(a^{k+1})|\} \right)^{k-1}. \end{align*}
Fixing $C_2 := 1$ and recursively taking $C_{k+1} := (k+1) C_k + k$, the constants $C_k$ must equal $k! - 1$ for each $k$, giving
\eqref{theinequality}. This proves the proposition.

\section{Proof of Theorem \ref{thm_main}}
To prove the theorem, we apply the pointwise inequality (\ref{real_inequality}) with $k=2r$. For simplicity, assume temporarily that the collection $J$ is finite. For a fixed $x \in X$, each sequence $a^i$ with $1 \leq i \leq r$ is chosen to be the sequence $\{f_j(x)\}_{j \in J}$; each sequence $a^i$ with $r+1 \leq i \leq 2r$ chosen to be the sequence $\{\overline{f_j}(x)\}_{j \in J}$. For each $x \in X$ we may then apply the real inequality (\ref{real_inequality}). By integrating the   pointwise inequality (\ref{real_inequality}) over $X$ we obtain
\begin{multline*}
  \| \sum_{j\in J} f_j \|^{2r}_{L^{2r}(X,d\mu)}
\leq ((2r)!-1) \int(B^{2r} + B^2A^{2r-2})  \\+ \Re( \sum_{\substack{j_1,\ldots,j_{2r}\in J \\ j_1,\ldots,j_{2r} \text{ distinct}}} \int f_{j_1} \cdots f_{j_r} \overline{f_{j_{r+1}} \cdots f_{j_{2r}}}).
\end{multline*}
The last term vanishes, by hypothesis. By the definition of $A,B$ we conclude that 
\beq\label{application}  \| \sum_{j\in J} f_j \|^{2r}_{L^{2r}(X,d\mu)}
\leq ((2r)!-1)\left( \int (\sum_{j\in J} |f_j|^2)^{r}
    + \int (\sum_{j\in J} |f_j|^2)|\sum_{j\in J} f_j|^{2r-2}\right).\eeq
The first term on the right-hand side is already of the form $ \| (\sum_{j\in J} |f_j|^2)^{1/2} \|^{2r}_{L^{2r}}$. 
To the second term on the right-hand side we apply H\"older's inequality, so that it is bounded above by 
\begin{multline*}
    (\int |\sum_{j\in J} f_j|^{2r})^{1 - 1/r} (\int (\sum_{j\in J} |f_j|^2)^{r})^{1/r}
        \\ \leq
 \frac{\ep^{r/(r-1)}}{r/(r-1)}\|\sum_{j\in J} f_j\|_{L^{2r}}^{2r} 
 + \frac{1}{\ep^{r}r} \| (\sum_{j\in J}|f_j|^2)^{1/2}\|_{L^{2r}}^{2r},
 \end{multline*} 
for any $\ep>0$. The last inequality follows from Young's inequality for products, which shows that for any conjugate exponents $1/p + 1/q=1$, and any real $a,b \geq 0$,
\[ a^{1/p}b^{1/q} \leq \ep^p\frac{ a}{p} + \ep^{-q}\frac{b}{q}, \quad \text{for any $\ep>0$.}\]
Here we notice that we have recovered a term of the form $\|\sum_{j\in J} f_j\|_{L^{2r}}^{2r}$ on the right hand side, and by choosing $\ep$ sufficiently small, we can subtract this term harmlessly from the left-hand side of (\ref{application}). Upon choosing $\ep$ so that 
\[ \frac{\ep^{r/(r-1)}}{r/(r-1)}((2r)!-1) \leq 1/2,\]
say, 
we conclude that 
\[  \frac{1}{2}\| \sum_{j\in J} f_j \|^{2r}_{L^{2r}(X,d\mu)}
\leq ((2r)!-1) \left( 1 + \frac{1}{\ep^r r} \right)\| (\sum_{j\in J} |f_j|^2)^{1/2}\|_{L^{2r}}^{2r}.\]
This proves the theorem for finite $J$. In particular we can observe that this proof allows
$C_r \leq 2^{1/2}((2r)!-1)^{1/2}
$.

When $J$ is infinite, \eqref{square} may be deduced from the finite case by an appropriate limiting argument. First observe that the pointwise values of the square function $(\sum_{j \in J} |f_j|^2)^{1/2}$ are independent of any choice of ordering of $J$ because the terms of the sum are nonnegative; we now assume it belongs to $L^{2r}$. Let $E$ and $E'$ be any finite subsets of $J$; then
\begin{align*} \| \sum_{j \in E} f_{j} - \sum_{j' \in E'} f_{j'} \|_{L^{2r}} & = \| \sum_{j \in E \setminus E'} f_{j} - \sum_{j' \in E' \setminus E} f_{j'} \|_{L^{2r}} \\ & \leq
\| \sum_{j \in E \setminus E'} f_{j_i}\|_{L^{2r}} + \| \sum_{j' \in E' \setminus E} f_{j'} \|_{L^{2r}} \\
& \leq C_r \| (\sum_{j \in E \setminus E'} |f_j|^2)^{\frac{1}{2}} \|_{L^{2r}} +  C_r \| (\sum_{j' \in E' \setminus E} |f_{j'}|^2)^{\frac{1}{2}} \|_{L^{2r}} \end{align*}
by applying the theorem for finite index sets.
If $j_1,j_2,\ldots$ and $j_1',j_2',\ldots$ are any enumerations of $J$, the Lebesgue Dominated Convergence Theorem guarantees that
$\| (\sum_{i = N+1}^\infty |f_{j_i}|^2)^{{1}/{2}} \|_{L^{2r}} \rightarrow 0$
as $N \rightarrow \infty$ and likewise $(\sum_{i' = N'+1}^\infty |f_{j'_{i'}}|^2)^{{1}/{2}}$ tends to zero in $L^{2r}$ as $N' \rightarrow \infty$. Taking $E = \{j_1,\ldots,j_N\}$ and $E' = \{j'_1,\ldots,j'_{N'}\}$ yields
\begin{align*} \| \sum_{i=1}^{N} f_{j_i} - \sum_{i'=1}^{N'} f_{j'_{i'}} \|_{L^{2r}} 
& \leq C_r \| (\sum_{j \in E \setminus E'} |f_j|^2)^{\frac{1}{2}} \|_{L^{2r}} +  C_r \| (\sum_{j' \in E' \setminus E} |f_{j'}|^2)^{\frac{1}{2}} \|_{L^{2r}} \\
& \leq C_r \| (\sum_{i=N+1}^\infty |f_{j_i}|^2)^{\frac{1}{2}} \|_{L^{2r}} +  C_r \| (\sum_{i' = N'+1}^\infty |f_{j'_{i'}}|^2)^{\frac{1}{2}} \|_{L^{2r}} \end{align*}
which implies that the partial sums are a Cauchy sequence (by taking the orderings $j_1,j_2,\ldots$ and $j'_1,j'_2,\ldots$ to coincide) and that all orderings of the series converge to the same limit in $L^{2r}$. Consequently
\begin{align*} \| \sum_{j\in J} f_j \|_{L^{2r}} = \lim_{N\rightarrow \infty} \| \sum_{i=1}^N f_{j_i} \|_{L^{2r}} \leq C_r \lim_{N\rightarrow \infty} \| (\sum_{i=1}^N |f_{j_i}|^2 )^{\frac{1}{2}}\|_{L^{2r}} = C_r \| (\sum_{j\in J} |f_j|^2)^\frac{1}{2}\|_{L^{2r}} \end{align*}
as desired.

\section{Construction of an example: proof of Theorem \ref{thm_example}}

To prove Theorem \ref{thm_example}, when $k=1$, the functions $f_i$ can be taken to be one of any number of piecewise-constant mutually orthogonal functions in $L^2(\R)$. So without loss of generality, we may assume $k \geq 2$. Let ${\mathcal I} \subset \N^{2k-2}$ consist of all tuples $(j_1,\ldots,j_{2k-2})$ of strictly-increasing indices: $1 \leq j_1 < j_2 < \cdots < j_{2k-2}$. Subdivide $\R$ into intervals of the form $[\ell,\ell+1)$ for each $\ell \in {\mathbb Z}$ and suppose that these intervals are indexed by elements of $\mathcal I$; in other words, assume that for each $(j_1,\ldots,j_{2k-2}) \in \mathcal I$, there is some interval $I_{j_1,\ldots,j_{2k-2}}$ having the form $[\ell,\ell+1)$ for some $\ell \in {\mathbb Z}$ such that no distinct elements of $\mathcal I$ are associated to the same such interval in $\R$ and every such interval $[\ell,\ell+1)$ is associated to exactly one element of $\mathcal I$.

For each $j \in \N$, let $r_j(t)$ be the $1$-periodic function on $\R$ which agrees with the $j$-th nonconstant Rademacher function on $[0,1)$ and let $g_j(t)$ be any function in $L^{2k}(\R)$ which is constant and strictly positive on each interval of the form $[\ell,\ell+1)$ for $\ell \in \Z$.
Now for each $i \in \N$,
define, for $t \in I_{j_1,\dots,j_{2k-2}}$,
\[ f_i(t) := \begin{cases} g_i(t) r_i(t) & i \not \in \{j_1,\ldots,j_{2k-2}\}, \\ g_i(t) & i \in \{j_1,\ldots,j_{2k-2}\}. \end{cases} \]
These functions $\{f_i\}_{i=1}^\infty$ will be shown to satisfy the conclusions of the theorem.

First consider \eqref{innerp0}. It suffices to show that for any $1 \leq i_1 < i_2 < \cdots < i_{2k}$ and any interval $[\ell,\ell+1)$,
\begin{equation} \int_{\ell}^{\ell+1} f_{i_1}(t) \cdots f_{i_{2k}}(t)\,dt = 0. \label{integralzero} \end{equation}
Each such interval $[\ell,\ell+1)$ is equal to $I_{j_1,\ldots,j_{2k-2}}$ for some indices $1 \leq j_1 < j_2 < \cdots < j_{2k-2}$. Because all the $i$'s are distinct, there must exist some $s \in \{1,\ldots,2k\}$ such that $i_s \neq j_{s'}$ for any $s' \in \{1,\ldots,2k-2\}$. This means that, on the interval $[\ell,\ell+1)$, $f_{i_1} \cdots f_{i_{2k}}$ is a constant multiple of some product of Rademacher functions, and in this product, the Rademacher function $r_{i_s}$ appears exactly once. Thus the orthogonality properties of the Rademacher functions imply that \eqref{integralzero} holds.

For \eqref{innerp1}, suppose that $i_1,\dots,i_{2k}$ are not distinct. On every interval $[\ell,\ell+1)$, $f_{i_1} \cdots f_{i_{2k}}$ is a strictly positive multiple of some product of Rademacher functions. As such, it is always the case that
\begin{equation} \int_{\ell}^{\ell+1} f_{i_1}(t) \cdots f_{i_{2k}}(t)\,dt \geq 0. \label{integralnotzero} \end{equation}
It suffices to show, then, that the integral \eqref{integralnotzero} is strictly positive for at least one $\ell$. To that end, suppose without loss of generality that $1 \leq i_1 \leq i_2 \leq \cdots \leq i_{2k}$. Let $j_1 < \cdots < j_m$ enumerate those natural numbers which appear in the tuple $(i_1,\ldots,i_{2k})$ an odd number of times. By assumption, there is at least one natural number appearing in this tuple two or more times, so the total number of natural numbers appearing in the list an odd number of times can be at most $2k-2$. (If there are any natural numbers at all appearing an even number of times, removing them from the list leaves at most $2k-2$ elements; and if all indices appear an odd number of times, at least one must appear 3 or more times, so removing two of its instances reduces the list length to $2k-2$ again while preserving the criterion by which $j_1,\ldots,j_m$ are chosen.) Extending the sequence of $j$'s as needed, we may assume that there exist $1 \leq j_1 < \cdots < j_{2k-2}$ such that every index $i$ appearing in the tuple $(i_1,\ldots,i_{2k})$ an odd number of times also appears in $(j_1,\ldots,j_{2k-2})$. Now consider the integral of $f_{i_1} (t)\cdots f_{i_{2k}}(t)$ on $I_{j_1,\ldots,j_{2k-2}}$. Any index $i$ appearing an odd number of times in $(i_1,\ldots,i_{2k})$ must belong to $\{j_1,\ldots,j_{2k-2}\}$ and consequently $f_i(t)$ is positive and constant on $I_{j_1,\ldots,j_{2k-2}}$ for each such $i$. If $i$ appears an even number of times in $(i_1,\ldots,i_{2k})$, then $f_i$ is a positive multiple of $r_i(t)$, and so raising $f_i$ to an even power makes it again strictly positive and constant. As a consequence, $f_{i_1} \cdots f_{i_{2k}}$ is strictly positive and constant on $I_{j_1,\ldots,j_{2k-2}}$ and therefore \eqref{integralnotzero} must be a strict inequality on $I_{j_1,\ldots,j_{2k-2}}$ as promised.

\section{Haar functions: Proposition \ref{prop_Haar}}
Recall the Haar functions $\Psi_I$ as defined in (\ref{Haar_dfn}). 
Observe that for $n$ dyadic intervals $I_1\subseteq \cdots \subseteq I_n\subsetneq I_{n+1}=\mathbb{R}$ we have
\[ \int \psi_{I_1} \cdots \psi_{I_n} \not= 0 \]
if and only if the number $m\in \{1,\dots,n\}$ defined by $I_1=\cdots=I_m\subsetneq I_{m+1}$ is even.
As a consequence, the family of Haar functions $\{\psi_I\}_{I\in\mathcal{D}}$, indexed by the set of dyadic intervals, is superorthogonal of Type IV:
if $I_1,\dots,I_{2r}$ are pairwise distinct, then $\int \psi_{I_1} \cdots \psi_{I_{2r}} = 0$ vanishes.
Haar functions are however not superorthogonal of Type II for $2r$-tuples with $r\ge 2$: if $I_1\subsetneq I_2\subsetneq I_3$ are dyadic intervals, then the index $I_3$ appears precisely once in the $2r$-tuple $(I_1,I_1,I_2,\dots,I_2,I_3)$, but
\[ \int \psi_{I_1}^2\; \psi_{I_2}^{2r-3}\; \psi_{I_3} \not= 0. \]
A similar argument shows that Haar functions are also not of Type III on $\R$. To be of Type III we would need to fix an ordering of the set $\mathcal{D}$, say a bijection $\iota: \mathbb{N}\to \mathcal{D}$ and set $\psi_{\iota(j)}=\psi_{I_j}$. According to this ordering, we would say the sequence $\{\psi_j\}_{j=1}^\infty$ is superorthogonal of Type III on $\R$ for some $r\ge 2$, in the sense that 
\beq\label{Haar_vanish} \int \psi_{I_{j_1}}\cdots \psi_{I_{j_{2r}}} =0
\eeq
whenever there is an index $j_i$ appearing in the product such that $j_i>j'$ for all $j' \neq j_i$, $j' \in \{j_1,\ldots,{j_{2r}}\}$. The obstacle is that no such ordering can exist.

Indeed, choose an arbitrary bijection $\iota: \mathbb{N}\to \mathcal{D}$ and set $I_j=\iota(j)$.
We claim that if $I_{j_1}\subsetneq I_{j_2}$ for two dyadic intervals and the Type III property (\ref{Haar_vanish}) holds, then necessarily $j_1>j_2$.
Indeed, suppose $I_{j_1}\subsetneq I_{j_2}$ are such that $j_1<j_2$. Choose a dyadic interval $I_{j_3}\supsetneq I_{j_2}$ and
let $i_{\mathrm{min}}=\min(j_2,j_3)$ and $i_{\mathrm{max}}=\max(j_2,j_3)>j_1$.
Then $i_{\mathrm{max}}$ is the largest index in the $2r$-tuple $(j_1, j_1, i_{\mathrm{min}}, \dots, i_{\mathrm{min}}, i_{\mathrm{max}})$, but
\[ \int \psi_{I_{j_1}}^2 \psi_{I_{i_{\mathrm{min}}}}^{2r-3} \psi_{I_{i_{\mathrm{max}}}} \not= 0, \]
because $I_{j_1}\subsetneq I_{j_2}\subsetneq I_{j_3}$.
This is a contradiction to the Type III property (\ref{Haar_vanish})   for this ordering. Thus, we must have $j_1>j_2$ whenever $I_{j_1}\subsetneq I_{j_2}$.
Such a bijection $\iota$ cannot exist, because the natural numbers are bounded from below (consider an ancestor of $I_1=\iota(1)$).

However, the subfamily of Haar functions associated with dyadic intervals contained in a given \emph{compact} interval is superorthogonal of Type III,  by indexing the sets according to $I_{(k,\ell)}=[2^k \ell, 2^k (\ell+1))$, and applying   lexicographic ordering of $(k,\ell),$ so that $(k,\ell) > (k',\ell')$ if $k>k'$ or $k=k'$ and $\ell>\ell'$.

\section{Further directions}

The  identification of superorthogonality as a formal tool for obtaining a square function estimate is quite new, and in the spirit of exploration we raise several questions.

Formalizing a hierarchy of Types is useful because it provides a variety  of conditions one can test for a given sequence of functions for which one wishes to prove a square function estimate. In particular, Type IV superorthogonality, and its unconditional implication of a square function estimate, promises to be useful  since one only needs to check (\ref{vanish'}) on a relatively small set, call it
\[ Z(\mathrm{IV})  = \{\text{$2r$-tuples with $j_1,\ldots,j_{2r}$ all distinct}\}.\]
Is there a Type of superorthogonality that unconditionally guarantees a square function estimate, with a   test set of index tuples that is strictly smaller than 
the set $Z(\mathrm{IV}) $?

We proved Theorem \ref{thm_main} with a constant $C_r \leq 2^{1/2} ((2r)!-1)^{1/2}$ for the square function estimate in $L^{2r}$, $r \geq 2$. Showing a sequence of functions is Type IV provides a direct route to proving a square function estimate by testing a condition on a relatively small  set of indices, but proving the sequence satisfies a more restrictive Type in the hierarchy can significantly improve the constant playing the role of $C_r$. For example, recently Hickman and Wright \cite{HicWri22x} and Hughes \cite{Hug22x} have  studied an extension operator associated to the moment curve in $n$-dimensional space, in the setting of  non-archimedean local fields. This type of study can be seen as a non-archimedean analogue of the    square function estimate obtained over $\R$ in the authors' work \cite{GGPRY20}. The latter work effectively obtained the square function estimate by proving that along a nondegenerate curve such as the moment curve, extension operators for sufficiently small, sufficiently separated intervals satisfy Type I* superorthogonality, so that a square function estimate holds with constant $C_r \leq (r!)^{1/2r}$ for each integer $r \leq n$.  Interestingly, without applying the full strength of their method,  Hickman and Wright prove that such extension operators satisfy a Type II condition (see \cite[Eqn(13)]{HicWri22x} and the two subsequent bullet points). By \cite[\S 3.1]{Pie20}, Type II functions satisfy a square function estimate, but the proof there only gives $C_r \leq r$. By pushing further to show a stronger Type I* property holds (via a non-archimedean version of the Phong-Stein-Sturm sublevel set decomposition, see \cite[Prop. 1.2]{HicWri22x}), Hickman and Wright's method allows the smaller choice $C_r \leq (r!)^{1/2r}$ for each integer $r \leq n$.   (Alternatively, Hughes verifies a Type I* property  using the Girard-Newton equations, and is able to take $C_r \leq (r^r)^{1/2r}=r^{1/2}$ for each integer $r \leq n$ \cite[Prop. 5]{Hug22x}.)
This leads to a question: can  the formal constant in the square function inequality implied by a Type of superorthogonality, such as $C_r$ in Theorem \ref{thm_main} for Type IV, be significantly improved in general?

 Pierce   also introduced a notion of quasi-superorthogonality \cite[\S 7]{Pie20}. Quasi-superorthogonality can be interpreted as  the condition that the integral (\ref{vanish'}) satisfies a  non-trivial upper bound (rather than vanishing) whenever the tuple of indices $j_1,\ldots,j_{2r}$ lies in a certain set. (It turns out that this notion cleanly characterizes an important phenomenon observed for multiplicative Dirichlet characters, and more generally trace functions, in analytic number theory.)
A different notion of being ``almost superorthogonal'' of a certain Type could be the requirement that (\ref{vanish'}) vanishes whenever the tuple of indices $j_1,\ldots,j_{2r}$ is ``close'' to the set $Z(\mathrm{Type})$. As a particularly simple example, being almost Type IV could impose that (\ref{vanish'}) vanishes for all $j_1,\ldots,j_{2r}$ such that  $\min_{i,i'} |j_i - j_{i'}| \geq c$. This simple case is reminiscent of   C\'ordoba's classic $L^4$ argument for Bochner-Riesz operators \cite[p. 507]{Cor79}. Or, Andreas Seeger has pointed out that both these notions can be seen at play in  \cite[Prop. 5.1]{ChrSee06}.  For some notions of ``almost,'' one can do a preliminary step that sparsifies or separates an original sequence $\{f_j\}_j$ that is ``almost'' a Type into a finite number of sequences, each of which is precisely that Type; this is demonstrated for Type I* in \cite[\S 6]{GGPRY20}. (See also Type III' in \cite[\S 4.5]{Pie20}.)  Are there   applications that essentially require interesting notions of ``almost'' or ``quasi'' superorthogonality?

\subsection*{Acknowledgements}
The authors thank J. Hickman and J. Wright for their encouragement and an interesting discussion at Oberwolfach in July 2022.
We thank AIM for funding our SQuaRE workshop. Gressman was partially supported by NSF DMS-1764143, DMS-2054602.  
Pierce  was partially supported by NSF CAREER DMS-1652173, DMS-2200470, a Sloan Research Fellowship, and a Joan and Joseph Birman Fellowship  during portions of this work, and thanks the Hausdorff Center for Mathematics for a productive visit as a Bonn Research Chair in 2022. Roos was partially supported by NSF DMS-2154356 and a grant from the Simons Foundation (ID 855692). Yung  was partially supported  by a Future Fellowship FT200100399 from the Australian Research Council.
%***************************************
\bibliographystyle{alpha}
\bibliography{Bib}
%***************************************

\end{document}